\documentclass{amsart}

\usepackage[pdftex]{graphicx} 
\usepackage[usenames,dvipsnames]{xcolor}

\usepackage[english]{babel} \usepackage[utf8]{inputenc}
\usepackage{tikz-cd}
\usepackage[centering]{geometry}
\usepackage{upgreek}
\usepackage{amsfonts}
\usepackage[font=small,labelfont=bf]{caption}
\usepackage{tikz}
\usepackage{setspace}
\usepackage{tensor}
\usepackage{lmodern}
\usepackage{verbatim}
\usepackage{amsmath}
\usepackage[all, cmtip]{xy}\usepackage{amssymb}
\usepackage{amsthm}
\usepackage{extpfeil}
\usepackage{bbm}
\usepackage{paralist}
\usepackage[colorlinks,citecolor=magenta,pagebackref=true,urlcolor=magenta,pdftex]{hyperref}
\usepackage{cancel}

\usepackage{soul}

\usepackage[T1]{fontenc}

\usepackage{nicefrac}
\usepackage{microtype}
\usepackage{mathrsfs}

\makeatletter
\newtheorem*{rep@theorem}{\rep@title}
\newcommand{\newreptheorem}[2]{%
\newenvironment{rep#1}[1]{%
 \def\rep@title{#2~\ref{##1}}%
 \begin{rep@theorem}}%
 {\end{rep@theorem}}}
 
\makeatother
\theoremstyle{plain}

\newtheorem*{thm*}{Theorem}

\newtheorem*{LBTP}{Lower Bound Theorem for polytopes}

\newtheorem{thm}{Theorem}[section]
\newtheorem{cor}[thm]{Corollary}
\newtheorem{lem}[thm]{Lemma}
\newtheorem*{lem*}{Lemma}

\newreptheorem{thm}{Theorem}
\newreptheorem{cor}{Corollary}

\theoremstyle{definition}

\newcommand*\longhookrightarrow{\ensuremath{\lhook\joinrel\relbar\joinrel\rightarrow}}

\definecolor{light-gray}{gray}{0.55}

\newcommand\hlightb[1]{\tikz[overlay, remember picture,baseline=-\the\dimexpr\fontdimen22\textfont2\relax]\node[rectangle,fill=blue!50,rounded corners,fill opacity = 0.2,draw,thick,text opacity =1] {$#1$};} 
\newcommand\hlightr[1]{\tikz[overlay, remember picture,baseline=-\the\dimexpr\fontdimen22\textfont2\relax]\node[rectangle,fill=red!50,rounded corners,fill opacity = 0.2,draw,thick,text opacity =1] {$#1$};}

\newcommand{\xqedhere}[2]{%
  \rlap{\hbox to#1{\hfil\llap{\ensuremath{#2}}}}}

\newcommand{\cm}[1]{}

\newcommand\mr[1]{\mathrm{#1}}

\newcommand{\R}{\mathbb{R}}

\DeclareMathOperator{\St}{st}

\title{Lefschetz and Lower Bound Theorems for Minkowski
sums}

\author{Karim Adiprasito}
\address{Einstein Institute for Mathematics, Hebrew University of Jerusalem, Jerusalem, Israel
}
\email{adiprasito@math.huji.ac.il}

 %  Combinatorial aspects of simplicial complexes

\begin{document}

\begin{abstract}This note provides a Lefschetz theorem for Minkowski sums of polytopes, and conclude lower bound theorems for Minkowski sums of polytopes. It is written as an appendix to arXiv:1405.7368, so notation and references follow that paper.
\end{abstract}

%{\color{magenta}
%Notation:
%\begin{compactitem}[$-$]
%    \item $A \subseteq B$: $A$ is a subset of $B$
%    %\item $A \subseteq B$: $A$ is a subset of $B$ and $A \neq B$
%    \item $A \subsetneq B$: $A$ is a subset of $B$ and $A \neq B$
%    \item $A \not\subseteq B$: $A$ is not a subset of $B$
%\end{compactitem}
%}

\maketitle

\newcommand\Neighborly{\mr{NP}}
\newcommand\Cyclic{\mr{Cyc}}
\newcommand\n{\mathbf{n}}

\vskip -10mm

\enlargethispage{10mm}

\setcounter{section}{7}

\section{Lower bounds for Minkowski sums}

Let us recall the following theorem due to Barnette. 

\begin{LBTP}
    For a simplicial $d$-dimensional polytope $P$ on $n$ vertices and $0 \le k < d$
    \[
        f_k(P) \ \ge \ f_k(\text{Stack}_d(n))
    \]
    where $\text{Stack}_d(n)$ is a $d$-dimensional stacked polytope on $n$ vertices. 
    Moreover, equality holds for all $k$ whenever it holds for some $k_0 \ge 1$, and the equality cases are attained by stacked polytopes only.
\end{LBTP}

This theorem was proven by Barnette \cite{Barnette}, with the equality case proven by Kalai \cite{KalaiRig}, and is substantially deeper 
than the Upper Bound Theorem as it relates to the Lefschetz property for toric varieties, see also \cite{KalaiRig, MN, AdiprasitoTC}.

In this section we will address more general lower bound problems for polytopes 
and polytopal complexes. 

{\it
For given $k < d$ and $n_1,n_2,\dots,n_m$, what is the \textbf{minimal} number of 
$k$-dimensional faces of the Minkowski sum $P_1 + P_2+\cdots+P_m$  for
polytopes $P_1,\dots,P_m \subseteq \R^d$ with vertex numbers $f_0(P_i) =
n_i$ for $i=1,\dots,m$?
}

To make this problem nontrivial, it is useful to restrict to "general position" Minkowski sums of "general position" polytopes, i.e., restrict to situations in which the Cayley complex is simplicial.   

For this, let us consider a family $P_{[m]} = (P_1,\dots,P_m)$ of simplicial polytopes in relatively general position in $\R^d$ with $\n = f_0(P_{[m]})$. For simplicity, we shall denote the Artinian reduction $\R[\Psi]/\Theta \R[\Psi]$ of a face ring or module $\R[\Psi]$ by $A(\Psi)$, and shall assume that $m\ge 2$, as the case $m=1$ is covered by the classical lower bound theorem for simplicial polytopes. We shall also assume that the Minkowski sum is of full dimension, as the question is trivial otherwise: the Minkowski sum is but a direct sum of simplices if the sum is of lower dimension.

\subsection{Lefschetz theorems for Minkowski sums}   
Let us start by examining the case of the minimal number of vertices in the Minkowski sum. That amounts to providing a lower bound for
$\widetilde{h}_{m}(\mr{T}_{[m]}^\circ)$ or equivalently, by Dehn-Sommerville,  ${h}_{d-1}(\mr{T}_{[m]})$ (with a correction term as in Dehn-Sommerville for nonpure collections). The key algebraic lemma is the following:

\begin{lem}
\[A(\mr{T}_{[m]})_k\ \xrightarrow{\ \cdot \ell^{d-2k} \ }\ A(\mr{T}_{[m]})_{d-k},\]
for $k\le \frac{d}{2}$, is an injection.
\end{lem}

Here, $\ell$ is an ample class on the Cayley polytope. 
\begin{proof} The lemma is a simple consequence of the Lefschetz theorem for simplicial polytopes, and the partition lemma \cite[Theorem 5.1]{AY}: We have an injection
\[A(\mr{T}_{[m]})_i\ \longhookrightarrow \ \bigoplus_{v \in \mr{T}_{[m]}} A(\St(v, \mr{T}_{[m]}) )_{i},\]
for $i\ <\ {d}$, where the direct sum is over vertices in $v$ in $\mr{T}_{[m]}$.

 It then follows that
\[\begin{tikzcd}
 A(\mr{T}_{[m]})_k\arrow{r}{\ \cdot  \ell^{d-2k}\ } \arrow[hook]{d}{} & A(\mr{T}_{[m]})_{d-k} \arrow[hook]{d}{} \\
\bigoplus_{v \in \mr{T}_{[m]}}   A(\St(v, \mr{T}_{[m]}) )_{k} \arrow{r}{\ \cdot  \ell^{d-2k}\ } & \bigoplus_{v \in \mr{T}_{[m]}} A(\St(v, \mr{T}_{[m]}) )_{d-k}
\end{tikzcd}
\]
The bottom maps are injective by the Lefschetz theorem for polytopes, as they are disks with boundary induced by cardinality $(m-1)$-faces, see \cite[Theorem 9.3]{AY}. It follows that the same holds for the top map.
\end{proof}
In particular, we have
\[\R(\mr{T}_{[m]})_1\ \xrightarrow{\ \cdot \ell^{d-2} \ }\ \R(\mr{T}_{[m]})_{d-1}.\]

Note finally that ${h}_{1}(\mr{T}_{[m]})$ is $| \n|- d-m+1$, where $|\n | $ is the sum over the entries of $\n$ . Hence we obtain a lower bound on $\widetilde{h}_{m}(\mr{T}_{[m]}^\circ)$ as $|\n|- d-m+1$. It follows by \cite{KalaiRig} that this is attained only if the Cayley polytope of the family is stacked, that is, it has a triangulation without interior $d+m-3$ faces. In particular, each summand is stacked as well. Let us note a particular corollary that seems nonobvious to prove directly.

\begin{cor}
A Minkowski sum of polytopes of positive dimension in $\R^d$ in general position has at least as many vertices as the Minkowski sum of the same number of segments in $\R^d$ in general position.
\end{cor}

\subsection{Lower bounds on higher-dimensional faces}

We can also prove similar results for higher face numbers. Again, we have that, modulo correction term in the nonpure case, 
\[\widetilde{h}_{k+m-1}(\mr{T}_{[m]}^\circ) \ =\ {h}_{d-k}(\mr{T}_{[m]}),\]
and a Lefschetz injection for $k< \frac{d}{2}$
\[\R(\mr{T}_{[m]})_k\ \xrightarrow{\ \cdot \ell^{d-2k} \ }\ \R(\mr{T}_{[m]})_{d-k}.\]
This map is an isomorphism if and only if the Cayley polytope is $k$-stacked, that is, it has a triangulation without interior faces of codimension $k+1$, see \cite{AdiprasitoTC, MN}.

\newcommand{\etalchar}[1]{$^{#1}$}
\def\cprime{$'$}
\providecommand{\bysame}{\leavevmode\hbox to3em{\hrulefill}\thinspace}
\providecommand{\MR}{\relax\ifhmode\unskip\space\fi MR }
\providecommand{\MRhref}[2]{%
  \href{http://www.ams.org/mathscinet-getitem?mr=#1}{#2}
}
\providecommand{\href}[2]{#2}

\end{document}